\documentclass[a4paper,12pt]{article}
\usepackage{graphicx}
\usepackage{amssymb}
\usepackage{amsmath, amscd, amsthm}
\usepackage{mathrsfs}
\usepackage[all]{xy}
\usepackage{amsfonts}
\usepackage{amsrefs}
\usepackage[usenames]{color}
\usepackage[dvipsnames]{xcolor}
\begin{document}

\theoremstyle{remark}
\newtheorem{Remark}{Remark}
\newtheorem*{Remarks}{Remarks}
\newtheorem*{Proof}{Proof}
\theoremstyle{definition}
\newtheorem{Definition}{Definition}
\newtheorem*{Example}{Example}
\newtheorem*{Examples}{Examples}

\theoremstyle{plain}
\newtheorem{Theorem}{Theorem}
\newtheorem{Theorema}{Theorem}
\newtheorem{Proposition}{Proposition}
\newtheorem{Corollary}{Corollary}
\newtheorem{Lemma}{Lemma}

\renewcommand{\labelenumi}{\roman{enumi})}
\newcommand{\vs}{\vspace{5mm}}
\renewcommand{\baselinestretch}{1.10}
\setlength{\parskip}{3mm}
\setlength{\parindent}{0em}
\small\normalsize

\newcommand{\C}{\mathbb C}
\newcommand{\MM}{\mathfrak M}
\newcommand{\N}{\mathbb N}
\newcommand{\M}{\mathbb M}
\newcommand{\D}{\mathfrak{Diffc}}

\begin{center}
{\Large\sc On the Mellin transform of a $\mathcal D$--module.
\footnote{\today}
}\\ $\ $ \\
{ Ricardo Garc\'{\i}a L\'opez
\footnote{
 Departament de Matem\`atiques i Inform\`atica.
 Universitat de Barcelona,
 Gran Via, 585.
 E-08007, Barcelona, Spain.
 e-mail: ricardogarcia@ub.edu}
 \footnote{Supported by the spanish Ministerio de Ciencia e Innovaci\'on research project PID2019-104844GB-I00}
}
\end{center}

\begin{abstract}
Given a holonomic $\C[z,z^{-1}]\langle \partial_z\rangle$-module $\M$, following \cite{LS} one can 
consider its Mellin transform,
which is a difference system on the affine line over $\C$. In this note we prove a stationary phase formula, which shows 
that its formal behavior at infinity 
is determined by the local germs defined by $\M$ at its singular points. 
\end{abstract}

{\bf 1. Introduction.}

The Fourier transform for $\mathcal D$-modules has been extensively studied, the most precise results available are in dimension one, that is,
for holonomic modules over $\mathbb C [x]\langle \partial_x \rangle$, see \cite{Mal}. In the analogous situation for $\ell$-adic sheaves, 
G. Laumon defined in \cite{Lau} so-called local Fourier transformations, which are related to the global $\ell$-adic Fourier-Deligne transform via a stationary phase formula. These local transformations allowed him to give a product formula for local constants, a construction of the Artin representation in equal characteristic and a simplification of Deligne's proof of the Weil conjecture. 

Having Laumon's work as a guideline, local Fourier transforms have been defined in the $\mathcal D$-module setting (\cite{BE}, \cite{GL}, \cite{DMR}, \cite{Ari}), where they also satisfy a stationary phase formula (see \cite{GL}), albeit only at the formal level. Beyond this, Stokes structures have to be considered and the study becomes much more complicated, see \cite{HS}, \cite{AHMS}, \cite{Mo}, \cite{AK}. 

In \cite{LS}, F. Loeser and C. Sabbah defined the Mellin transform of a $\mathcal D$-module on an algebraic torus (see also \cite{Lau2}), and they used it to prove a product formula for the determinant of the Aomoto complex (\cite{LS}*{Th\'eor\`eme 2.3.1}).

In \cite{GS},  A. Graham-Squire defined local Mellin transforms for formal germs of meromorphic connections in one variable, and computed them explicitly. They might be regarded as local analogues of the global Mellin transform of Loeser and Sabbah, and in this note we prove that a stationary phase formula holds also in this case.

More precisely, if $\mathfrak H$ denotes the category of holonomic $\C[[x]]\langle\partial_x\rangle$-modules, $\mathfrak H'$ the category of formal connections and 
$\D$ the category of difference $\C((\theta))$-modules, then we define local Mellin transform functors
\begin{eqnarray*}
  \mathfrak M^{(s,\infty)}(\bullet): \mathfrak H &\longrightarrow & \D \ ,\ s\in\C\smallsetminus\{0\} \\
   \mathfrak M^{(\star,\infty)}(\bullet): \mathfrak H' &\longrightarrow & \D \ ,\ \star\in\{0,\infty\}\,.
\end{eqnarray*}
The definition we give is microlocal, in the spirit of \cite{GL}, and extends that in \cite{GS}, in the sense that it allows to remove the assumptions made in loc. cit. about slopes and non-existence of horizontal sections.\footnote{But then, while in \cite{GS} the local Mellin transforms were equivalences of categories, here they are not.}

Let $\M$ be a holonomic $\C[z,z^{-1}]\langle z\partial_z\rangle$-module with singular set $S(\M)\subset\mathbb{C}\smallsetminus\{0\}$. For $\star\in\C\cup\{\infty\}$, denote by $\M_{\star}$ the local $\mathcal D$-module germ defined by $\M$ at $\star$, denote $\MM(\M)$ the global Mellin transform
of $\M$. 
This is a difference module on the affine line with coordinate $\eta$, and we can consider
 its germ at infinity $\MM(\M)_{\infty}$, which is a $\C((\theta))$-difference module, where $\theta=\eta^{-1}$. Our main result is:

{\bf Theorem: }{\it 
There is a functorial isomorphism of $\C((\theta))$-difference modules
\[
\MM(\M)_{\infty} \stackrel{\sim}\longrightarrow  \!\!\! \!\!\!\!\bigoplus_{\star\in S(\mathbb M)\cup\{0,\infty\}}  \!\!\! \!\!\!\!\MM^{(\star, \infty)}(\M_{\star}).
\] }

That is, the local (formal) behavior at infinity of the global 
Mellin transform of  $\M$ is determined by the (formal) germs 
defined by $\M$ at zero, infinity and at its singular points, and no global information is required.

I thank Claude Sabbah for his careful reading of a previous version of this note and for pointing up some mistakes and inaccuracies.

Through this note, $\C$ will denote the field of complex numbers. If $x$ is a coordinate, we denote 
$K_x=\C[[x]][x^{-1}]=\C((x))$ the field of Laurent series with complex coefficients on the variable $x$ and we denote $\C[x]\langle\partial_x\rangle$ the Weyl algebra of differential operators with polynomial coefficients. Unless otherwise specified, by a module over a non-commutative ring we mean a left module. 

\vs

{\bf 2. Differential and difference modules.} 

We recall a few well-known notions and results from the local theory of $\mathcal D$-modules and difference modules, 
we refer to \cite{Mal}, \cite{Sab} and \cite{Pra} for more details and proofs:

If $R=\C[x]\langle\partial_x\rangle\,,  R=\C[[x]]\langle \partial_x\rangle\,, R=  \C[x,x^{-1}]\langle\partial_x\rangle\text{ or } R= K_x\langle \partial_x\rangle$, an  $R$-module $M$ is holonomic if there is a non-zero left ideal $I\subset R$ such that $M$ is isomorphic to 
$R/I$ as a $R$-module. Holonomic $\C[[x]]\langle \partial_x\rangle$-modules $M$ such that 
$M=M[x^{-1}]$ are finitely dimensional over $K_x$, and they 
will be called formal connections.\footnote{The equivalence of this definition with the usual one is shown for instance in \cite{Sab}*{Theorem 4.3.2}.} For any holonomic $M$,
the localization $M[x^{-1}]$ is a formal connection.


Using the cyclic vector lemma, one attaches to a formal connection $M$ its 
Newton polygon (see for example 
\cite{Mal}*{Chapter III}), the slopes of its non-vertical sides are called the  formal slopes of $M$.  
One has a  functorial decomposition
\[
 M =\bigoplus_{\lambda} M^{\lambda}\,,
\]
where $\lambda$ runs over the set of slopes of $M$ and $M^{\lambda}$ is a formal connection which has only slope $\lambda$.

If $M$ is a formal connection, we denote by $\mathrm{irr}(M)$ its irregularity, defined 
as the height of the Newton polygon of $M$
(\cite{Mal}*{Chapitre IV, (4.5)}), and by $\mu(M)$ the dimension of its space 
of vanishing cycles ([loc. cit., \S 4]). By [loc. cit., Chapitre IV, Corollaire 4.10], we have   $\mu(M)= \dim(M)+\mathrm{irr}(M)$.\footnote{For the purposes of this note, this equality can be taken as a definition of $\mu(M)$.} For a holonomic $\C[[x]]\langle \partial_x\rangle$-module $M$,  both $\mathrm{irr}(M)$ and $\mu(M)$  are defined as those of $M[x^{-1}]$.

\begin{Definition} Let $\theta$ be a coordinate, denote $\phi: K_{\theta} 
\longrightarrow K_{\theta}$ the automorphism 
given by $\phi(a(\theta))= a(\theta/1+\theta)$. A difference module $(V,\Phi)$ is a 
finite-dimensional $K_{\theta}$-vector space $V$ endowed with a $\C$-linear invertible 
operator  $\Phi: V \longrightarrow V$ 
such that, for all $f\in K_{\theta}$ and $v\in V$ one has
\[
\Phi(f(\theta)\cdot v)= f(\phi(\theta))\cdot \Phi(v).
\]
Taking as morphisms those $\C$-linear maps which commute with the difference operators, difference modules over $K_{\theta}$ form an abelian category.
\end{Definition}

We briefly recall the
construction of the Newton polygon attached to a difference operator (see e.g. \cite{Pra}): 
Let $K_{\theta}\langle \Phi\rangle$ denote the skew-polynomial ring determined 
by the relations 
$\Phi\cdot f=\phi(f)\cdot \Phi$ for $f\in K_{\theta}$. 
With respect to the degree function, $K_{\theta}\langle \Phi\rangle$ is an 
euclidean ring and every finitely 
generated $K_{\theta}\langle \Phi\rangle$-module is a direct sum 
of cyclic modules.
The datum of a difference module is equivalent to that of a $K_{\theta}\langle \Phi\rangle$-module, of finite dimension as a $K_{\theta}$-vector space and such that the action of $\Phi$ is invertible,\footnote{We will always assume invertibility of $\Phi$.} or to the datum of a $K_{\theta}\langle \Phi, \Phi^{-1}\rangle$-module, finitely dimensional over $K_{\theta}$.

Given $P=\sum_{i=0}^m a_i\Phi^i\in K_{\theta}\langle \Phi\rangle$, the Newton 
polygon $\mathcal N(P)$ of $P$ is 
the convex envelope in $\mathbb R^2$ of the union of the half--lines $\{(x, 
y)\in \mathbb R^2 \, \mid x=i\ , \, y\geqslant v(a_i) \}$,
where $v:K_{\theta} \longrightarrow \mathbb Z\cup\{\infty\}$ is the $\theta$-valuation,
given by 
$v\left(\sum_j \alpha_j\theta^j\right)=\min\{j\mid \alpha_j\neq 0\}$, 
$v(0)=\infty$.\footnote{In \cite{Pra}, the condition defining the half-lines is 
$y\leqslant v(a_i)$ which, in view of the claimed properties of slopes, seems to be a misprint. Notice also that Praagman's polygon is not identical to the one considered in \cite{Du}, they differ by a reflection.}

It is proved in \cite{Pra}*{pg. 257, Remark 3}
that, up to a vertical translation corresponding to multiplication by a power of $\theta$, 
the polygon $\mathcal N(P)$ depends only on the difference module $D_P=K_{\theta}\langle \Phi\rangle/ K_{\theta}\langle \Phi\rangle \cdot P$. 
In particular, it follows easily from the definitions that the width of $\mathcal N(P)$ 
coincides with the dimension of $D_P$ as a 
$K_{\theta}$-vector space.
\begin{center}
\includegraphics[height=4cm]{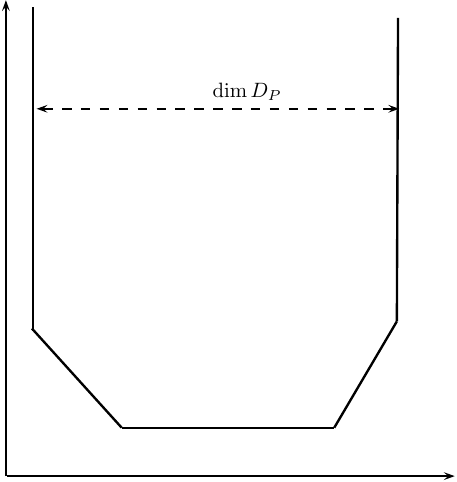}
\end{center}
In the sequel, the polygon $\mathcal N(P)$ will be always considered up to a vertical translation, the slopes of its non-vertical sides will be called the slopes of $D_P$.\footnote{ In \cite{GS} and \cite{Pra} the terminology differs. The orders considered by Graham-Squire are minus the slopes in Praagman's article.} In fact, if $V$ is a difference module over $K_{\theta}$, a version of the cyclic vector lemma allows to attach to $V$ an operator $P\in K_{\theta}\langle \Phi\rangle$ such that $V\cong D_P$ as difference modules. 

Given $q\geqslant 1$, set $L_q=\C((\theta^{1/q}))$. The automorphism $\phi(\theta)=\frac{\theta}{1+\theta}$ of $K_{\theta}$  
has a unique extension $\phi_q$ to $\C((\theta^{1/q}))$ (\cite{Pra}*{\S 1}), then one can define difference modules over $L_q$ as done for $K_{\theta}$, the definition of the Newton polygon of an operator extends as well.
Given $g\in L_q\smallsetminus \{0\}$,  we denote $D_{g,q}$ the $L_q$-difference module $(L_q, g\,\phi_q)$. Given $m\geqslant 0$, put $T_m=(K_{\theta}^ m, (Id+\theta N_m)\phi)$, where $N_m$ is the nilpotent 
Jordan block of size $m$. These are unipotent objects in the category of difference modules over $K_{\theta}$. The classification theorem for formal  difference modules is the following (see \cite{Pra}*{Theorem 8}, \cite{Du}*{Theorem 3.3}):

\begin{Theorem} \label{clas}
Let $V$ be a difference module over $K_{\theta}$ of dimension $m$. Then, there is a finite cyclic extension $K_{\theta}\subset L_q$ and an isomorphism of $L_q$-difference modules 
\[
 V\otimes L_q \cong\bigoplus_{i\in I} (D_{g_{i},q_{i}} \otimes_{K_{\theta}} T_{m_i})\otimes_{L_{q_i}}L_q\,,
\]
where $I$ is a finite set, $m_i, q_i>0$ are positive integers, $q_i\mid q$, the $D_{g_{i},q_{i}}$ are simple difference modules and $\sum m_i=m$. Also, $g_i\in \C((\theta^{1/q_i}))$ are of the form 
$g_i=\sum_{h=0}^{q_i}a_{i,h}\, \theta^{\lambda_i+\frac{h}{q_i}}$,
where $\lambda_i\in (1/q_i)\mathbb Z$ is the only slope of $D_{g_{i},q_{i}}$ and $a_{i,0}\in\C\smallsetminus \{0\}$. In this decomposition, the rational numbers $\{\lambda_i\}_{i\in I}$ are the formal slopes of $V$, the integers $m_i,q_i$ are 
uniquely determined, and the $g_i\in \C((\theta^{1/q_i}))$ are uniquely determined up to addition of an integer multiple of $a_{i,0}/q_i$ to $a_{i,q_i}\in \C$.

\end{Theorem}

It follows from the theorem that we have:
\begin{Corollary} \label{coslo}
If $V,W$ are difference modules with no common slope, every morphism  of difference modules $V \longrightarrow W$ is zero. 
\end{Corollary}
To prove a stationary phase formula for the Mellin transform, we will need some more information on the formal structure of a difference module than the one provided by formal slopes.\footnote{This is contrast with the situation 
for the stationary phase formula for the Fourier transform. In that case, formal slopes are enough, this difference is ultimately due to the different behavior of slopes with respect to tensor product in the differential and in the difference case.}

Let $P=\sum_{i\geqslant 0} a_{i}(\theta)\Phi^i\in K_{\theta}\langle \Phi\rangle$ with $a_i(\theta)=\sum_{i\geqslant 0}a_{ij}\theta^j\in K_{\theta}$, and assume $\mathcal N(P)$ has a 
horizontal side $\sigma$. Let $i_0<\dots < i_r$ be those indexes such that $a_{i_{\ell}}\Phi^{i_{\ell}}$ corresponds to a point on $\sigma$. For $0\leqslant \ell \leqslant r$, put $j_\ell=\min\{ j\in\mathbb N\,\mid\, a_{i_{\ell}\,j}\neq 0\}$ and consider the polynomial
\[
 p_{\sigma}(t)=\sum_{\ell=0}^r a_{i_\ell, j_\ell} t^{i_{\ell}-i_0}\in \C[t].
\]
 Then, the roots of $p_{\sigma}$ are formal invariants
of $D_P$, see e.g. \cite{BaCh}*{section 2.3}.
\begin{Definition}
 Let $V$ be a difference module over $K_{\theta}$, choose $P\in K_{\theta}\langle \Phi\rangle$ with $V\cong D_P$ as difference modules.
 We define a finite set $Hor(V)$ of horizontal zeros as follows: If the Newton polygon $\mathcal N(P)$ has no horizontal side, we put $Hor(V)=\emptyset$. If it has a horizontal side $\sigma_{hor}$, then $Hor(V)$ is the set of zeros  $p_{\sigma_{hor}}$, with multiplicities. This definition extends to the case $P\in L_q\langle \Phi\rangle$,  notice that for a $K_{\theta}$-difference module $V$ we have $Hor(V)=Hor(V\otimes L_q)$.
\end{Definition}

\begin{Lemma} \label{hor}
Let $V,W$ be difference modules over $K_{\theta}$, both of them purely of slope zero, with no common horizontal zero. 
Then every morphism $V\longrightarrow W$ is the zero map.
\end{Lemma}
\begin{Proof}
 Taking an extension $K_{\theta}\subset L_q$ we can assume both $V$ and $W$ decompose as in Theorem 1. 
Given a summand $D_{g,q}\otimes T_m$, a computation as in \cite{Pra} shows that 
\[
 D_{g,q}\otimes T_m \cong \frac{L_q\langle \Phi\rangle}{(\Phi-g)^m}
\]
where $g=a_0+a_1\theta^{1/q}+\dots+a_{q}\theta$ and $a_i\in\C$ for $0\leqslant i\leqslant q$. A direct calculation shows that the only horizontal zero of this difference module is $a_0\in\C$. Then, by the classification theorem, the lemma follows. $\Box$
\end{Proof}

\begin{Remark} \label{r1}
 It results from Theorem \ref{clas},  Corollary \ref{coslo} and Lemma \ref{hor} that if 
 \[
  0 \longrightarrow V' \longrightarrow V \longrightarrow V'' \longrightarrow 0
 \]
is an exact sequence of difference modules, then 
\begin{eqnarray}
 \{\text{Slopes of $V$}\} &= &\{\text{Slopes of $V'$}\}\cup  \{\text{Slopes of $V''$}\} \label{1} \\
 Hor(V) &=& Hor(V') \cup Hor(V''). \label{2}
\end{eqnarray}
The second equality is also valid when we consider the corresponding multiplicities.
In fact, if we define the multiplicity of a slope as the length of the vertical projection of the corresponding side onto 
the horizontal axis, then the first equality is also valid when we take multiplicities into account.

\end{Remark}

For later use, we recall a few definitions and results concerning differential and difference modules over tori and affine lines.
We recall that  $\C[z,z^{-1}]\langle \partial_z\rangle$ denotes the localized Weyl algebra
(where $[\partial_z,z]=1$) and $\C [\eta]\langle \Phi, \Phi^{-1}\rangle$ the 
algebra of 
invertible difference operators on the affine $\C$-line
(where $[\Phi,\eta]=\Phi$).

\begin{enumerate}
 \item If $\M$ is a holonomic $\C[z,z^{-1}]\langle 
\partial_z\rangle$-module and $s\in \C\smallsetminus\{0\}$, we set $z_s=z-s$ and 
we denote $\M_s$ the $\C[[z_s]]\langle \partial_{z_s}\rangle$-module $\C[[z_s]]\langle \partial_{z_s}\rangle\otimes_{\C[z,z^{-1}]\langle 
\partial_z\rangle}\M$.
We denote $\mathbb M_{0}$ the $K_{z}\langle \partial_z\rangle$-module $K_z\langle \partial_z\rangle\otimes_{\C[z,z^{-1}]\langle \partial_z\rangle}\M$ and, if $y=z^{-1}$, we denote $\mathbb M_{\infty}$ the $K_y\langle \partial_y\rangle$-module
$K_{y}\langle \partial_y\rangle\otimes_{\C[z,z^{-1}]\langle \partial_z\rangle}\M$, where   $\partial_y(1\otimes m)= 1\otimes (-z^2\partial_z )m$. Both $\M_0$ and $\M_{\infty}$ are formal connections.

The finite set of points $s\in \C\smallsetminus\{0\}$ such that $\C[[z_s]]\otimes_{\C[z]}\M$ is not free 
of finite type over 
$\C[[z_s]]$ will be denoted $S(\M)$ (the singular set of $\M$, see e.g. \cite{Sab}*{III. Proposition 1.1.5}).

\item The {\it global} Newton polygon attached by J. P. Ramis and B. Malgrange 
to  an operator 
$P\in \C[z,z^{-1}]\langle z\partial_z\rangle$ is 
defined as follows:\footnote{The definition in \cite{Mal}*{V.1} looks slightly 
different to the one given here, but both give the same polygon. In fact, in loc. cit. the case considered is that of a module over $\C[z]\langle \partial_z \rangle$. Since we have inverted $z$, we consider the global Newton polygon only up to horizontal translation.} Write $P$ as a finite sum
\[
 P=\sum_r \alpha_{r}(z\partial_z)z^r \text{ where } \alpha_{r}\in\C[X] \text{ for all $r\in\mathbb Z$}
\]
and, for each $\alpha_{r}\neq 0$, consider the half--line $\{(u,v)\in 
\mathbb 
R^2 \mid u\leqslant \deg{\alpha_r}, v=r\}$. 
The Newton polygon $\mathcal N(P)$ of $P$ is the convex envelope of 
these half--lines. It depends only on the quotient module $\mathbb M= \C[z,z^{-1}]\langle z\partial_z\rangle/ \C[z,z^{-1}]\langle 
z\partial_z\rangle P$ and, in particular, it follows from the definitions (see \cite{Mal}*{V.1}) that
the height $h(P)$ of $\mathcal N(P)$ equals
\[
 \mathrm{irr}(\mathbb M_0) + \mathrm{irr}(\mathbb M_{\infty}) + \sum_{s\in S(\mathbb 
M)}\mu(\mathbb M_s).
\]
\begin{center}
\includegraphics[height=4cm]{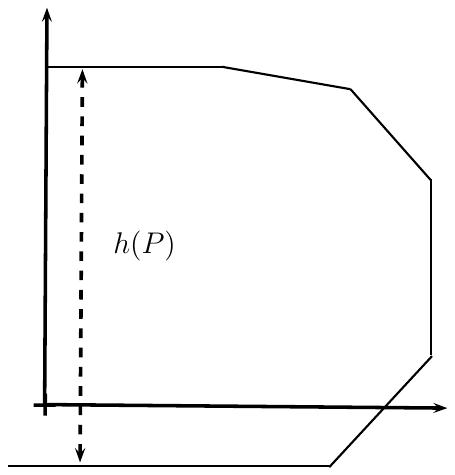}
\end{center}
As shown in loc. cit., in fact a Newton polygon can be attached to any holonomic $\C[z,z^{-1}]\langle z\partial_z\rangle$-module, and if 
\[
 0 \longrightarrow \M' \longrightarrow \M \longrightarrow \M'' \longrightarrow 0
\]
is an exact sequence, then 
\begin{equation} \label{3}
 \{\text{Slopes of $\M$}\} = \{\text{Slopes of $\M'$}\}\cup  \{\text{Slopes of $\M''$}\}, 
\end{equation}
also when we consider the slopes with their multiplicity, defined as the length of the horizontal projection of the corresponding sides onto the vertical axis .

\item If $\mathfrak N$ is a  $\C [\eta]\langle \Phi \rangle$-module, 
we define its germ at infinity as the $K_{\theta}$-vector space $\mathfrak N_{\infty}= K_{\theta}\otimes_{\C[\eta]} \mathfrak
N$, 
where $\eta \mapsto \theta^{-1}$, endowed with the difference operator given by
\[
 a(\theta)  \otimes n\longmapsto a\left(\frac{\theta}{1+\theta}\right)\otimes \Phi\cdot n  .
\]
Equivalently, $\mathfrak N_{\infty}$ is the $K_{\theta}\langle \Phi \rangle$-module  $K_{\theta}\langle \Phi \rangle \otimes_{\C [\eta]\langle \Phi \rangle} \mathfrak N$, or the $K_{\theta}\langle \Phi, \Phi^{-1}\rangle$-module $K_{\theta}\langle \Phi, \Phi^{-1}\rangle \otimes_{\C [\eta]\langle \Phi \rangle} \mathfrak N$.
\end{enumerate}

{\it Global Mellin transform:} Denote $\MM: 
\C[z,z^{-1}]\langle z\partial_z\rangle \longrightarrow 
\C [\eta]\langle \Phi, \Phi^{-1}\rangle$ the morphism of $\C$-algebras 
defined by $\MM (z\partial_z)= -\eta$, $\MM (z)= \Phi$. Following \cite{LS}, if 
$\M$ is a $\C[z,z^{-1}]\langle z\partial_z\rangle$-module, 
its Mellin transform is defined as the $\C [\eta]\langle \Phi, \Phi^{-1}\rangle$-module 
\[
 \MM(\M)=\C [\eta]\langle 
\Phi, \Phi^{-1}\rangle  \otimes_{\C[z,z^{-1}]\langle z\partial_z\rangle}\M.
\]
\begin{Remark}
i) If $\M$ is a $\C[z,z^{-1}]\langle z\partial_z\rangle$-module and $\N\subset \M$ is a $\C[z]\langle z\partial_z\rangle$-submodule which generates $\M$ over $\C[z,z^{-1}]\langle z\partial_z\rangle$, then we have
\begin{eqnarray*}
  \MM(\M)&\cong &\C [\eta]\langle 
\Phi, \Phi^{-1}\rangle \otimes_{\C[z,z^{-1}]\langle z\partial_z\rangle}  \C[z,z^{-1}]\langle z\partial_z\rangle\otimes_{\C[z]\langle z\partial_z\rangle}\N\\ &\cong&\C [\eta]\langle 
\Phi, \Phi^{-1}\rangle \otimes_{\C[z]\langle z\partial_z\rangle}\N \,,
\end{eqnarray*}
and similarly for a  $\C[z^{-1}]\langle z^{-1}\partial_{z^{-1}}\rangle$-submodule which generates $\M$ over \\ $\C[z,z^{-1}]\langle z\partial_z\rangle$.

ii) If $\M$ is holonomic, then $\MM(\M)_{\infty}$ is a quotient of $K_{\theta}\langle \Phi, \Phi^{-1}\rangle$ by a non-zero ideal, and so it is finite dimensional over $K_{\theta}$.

\end{Remark}

\begin{Lemma} \label{pol}
If $\M$ is a holonomic $\C[z,z^{-1}]\langle z\partial_z\rangle$-module, the Newton polygon of $ \MM(\M)_{\infty}$ is the polygon obtained from the Newton polygon of $\M$
by applying a rotation of ninety degrees in the clockwise direction.
\end{Lemma}
\begin{Proof}
 If $\M=\C[z,z^{-1}]\langle z\partial_z\rangle/ \C[z,z^{-1}]\langle z\partial_z\rangle\cdot P$, the claim is a consequence of the definitions, and also for modules with punctual support, which are a direct sum of modules of type $\C[z,z^{-1}]\langle z\partial_z\rangle/ \C[z,z^{-1}]\langle z\partial_z\rangle\cdot (z-s)^m$, with $s\in\C\smallsetminus\{0\}$ and $m\geqslant 1$. In general, there exists an exact sequence of $\C[z,z^{-1}]\langle z\partial_z\rangle$-modules
\[
 0 \longrightarrow \mathbb K \longrightarrow \M'\longrightarrow \M\longrightarrow 0\,,
\]
where $\mathbb K$ has punctual support (and so its Newton polygon consists only of a vertical side), and 
$\mathbb M'$ is a quotient by the left ideal generated by a single operator. By (\ref{3}) and the remark which follows it, the Newton polygon of $\M'$ coincides with that of $\M$, except for the fact that the vertical side of $\M'$ has as length the sum of the lengths of the vertical sides of $\M$ and of $\mathbb K$.  By exactness of the global Mellin transform and localization, we get an exact sequence
\[
 0 \longrightarrow \MM(\mathbb K)_{\infty} \longrightarrow \MM(\M')_{\infty}\longrightarrow \MM(\M^0_N)_{\infty}\longrightarrow 0.
\]
Then, the claim follows from Remark \ref{r1}. $\Box$
\end{Proof}

{\it Canonical good filtrations:} Let $R$ denote the ring $\C[z,z^{-1}]\langle z\partial_z \rangle$. Consider in 
$R$ the increasing filtration, indexed by $\mathbb Z$, given by 
\[
 V_kR=\left\{ \sum_i a_i(z,z^{-1})\partial_z^i \in R\, \mid \max_i\{ i-\operatorname{ord}_z(a_i)\}\leqslant k\right\}\,,
\]
and denote 
\[
  \Sigma=\{\alpha\in\C\,\mid \, -1\leqslant \operatorname{Re}\alpha\leqslant 0, \operatorname{Im}\alpha \geqslant 0 \text{ if } \operatorname{Re}\alpha=-1\ , \operatorname{Im}\alpha < 0 \text{ if } \operatorname{Re}\alpha= 0 \}.
 \]
If $\M$ is a holonomic $R$-module, there exists a unique increasing, $\mathbb Z$-indexed filtration $V_{\bullet}\M$ of $\M$ which is good for $V_{\bullet}R$ and such that the roots of its Bernstein polynomial 
are contained in $\Sigma$ (see \cite{Sab}*{Proposition 6.1.2}). This will be called the canonical good $z$-filtration of $\M$.\footnote{As in loc. cit., this choice of $\Sigma$ will play no special role in the sequel, it only matters that $\Sigma$ is a fundamental domain for the action of $\mathbb Z$ on $\C$ given by the translations $z\mapsto z+k, k\in\mathbb Z$.}

We can also consider the increasing filtration indexed by $\mathbb Z$ given by 
\[
\widetilde{V_kR}=\left\{ \sum_i a_i(z,z^{-1})\partial_z^i \in R\, \mid \max_i\{ i-\operatorname{ord}_{z^{-1}}(a_i)\}\leqslant k\right\}\,,
\]
and, in the same way, we obtain a filtration $\widetilde{V_{\bullet}\M}$ on $\M$ that we call its canonical good $z^{-1}$-filtration. 

One can also define a canonical good filtration for holonomic modules over $\C[[z]]\langle \partial_z\rangle$. 
If $\M$ is a holonomic $R$-module, then it is not difficult to see that we have
\[
  V_{\bullet}\M_0  =  \C[[z]]\otimes_{\C[z]}V_{\bullet}\M \ \ \text{ and } \ \
   V_{\bullet}\M_\infty  =  \C[[z]]\otimes_{\C[z]}\widetilde {V_{\bullet}\M}\,.
\]

%

{\bf 3. Microdifference operators and local Mellin transforms.}

In this section we will consider several rings, which differ by the product considered in each case. Set
\[
\mathfrak M = \left\{ \,\sum_{i\geqslant r} a_i(u) \eta^{-i}\ \mid r\in\mathbb Z \text{ and } a_i(u)\in\C[[u]] \,\right\}.
\]
The order of an operator $P\in\MM$  is  $\operatorname{ord}_{\eta}(P)=\max\{\,i \!\mid\! a_i\neq 0\,\}$. Let $\partial_{\eta}=d/d\eta$ denote the formal derivative with respect to $\eta$, let
$\delta$ be a $\C$--derivation of $\C[[u]]$. They can be formally extended  to maps defined in $\MM$,   
and we consider in $\MM$ 
the multiplication
\[
 P\circ_{\delta} Q=\sum_{\alpha\geqslant 0} \frac{1}{\alpha !}\ \partial_{\eta}^{\alpha}P\cdot \delta^{\alpha}Q\,,
\]
where, on the right hand side, it is understood that the product is first performed in  in $\C[[u]]((\eta^{-1}))$, and then the result is reordered so as to obtain an element of $\MM$.
It is easy to see that the set $\MM$, endowed with the obvious addition and the 
multiplication $\circ_{\delta}$, is a non-commutative unitary ring. 
We will only consider the derivations
$\delta_s=-(u+s)\partial_u$, where $s\in \C$, and $\delta_{\infty}=u\partial_u$. 

\begin{Definition} 
\begin{enumerate}
 \item For $s\in \C$, we denote by $\mathfrak M^{(s,\infty)}$ the ring which is $\MM$ as a set and where the product is
$\circ_{\delta_s}$. We regard it also as a $K_{\theta}$-algebra via the $\C$-morphism $K_{\theta}\longrightarrow \mathfrak M^{(s,\infty)}$ defined by $\theta \longmapsto -\eta^ {-1}$.
\item We denote by $\mathfrak M^{(\infty,\infty)}$ the set  
\[
 \left\{ \,\sum_{i\geqslant r} a_i(u^{-1}) \eta^{-i}\ \mid r\in\mathbb Z \text{ and } a_i(u^{-1})\in\C[[u^{-1}]] \,\right\}
\]
endowed with the product $\circ_{\delta_{\infty}}$ defined as above, where the derivation is $\delta_{\infty}=-u^{-1}\partial_{u^{-1}}(=u\partial_u)$
. We regard it also as a $K_{\theta}$-algebra via the morphism of $\C$-algebras
$K_{\theta} \longrightarrow \mathfrak M^{(\infty,\infty)}$   defined by $\theta \longmapsto \eta^ {-1}$.
\end{enumerate}
\end{Definition}

The proof of the following theorem is analogous to the one for the usual formal microdifferential operators
(see \cite{Du2}*{Th\'eor\`eme 1.c.2}):

\begin{Theorem} (division) For $\star\in\C\cup\{\infty \}$, set $u_{\star}=u+s$ if $\star=s\in\C$, $u_{\infty}=u^{-1}$. Let $P=\sum_{i\leqslant d} a_i(u_{\star})\eta^i\in \mathfrak M^{(\star,\infty)}$ be an operator of order $d$. Assume $a_d(u_{\star})\in\C[[u_{\star}]]$ has $u_{\star}$-adic valuation $p\in\mathbb N$.
Then, for all $S\in \mathfrak M^{(\star,\infty)}$ there are unique $Q,R\in \mathfrak M^{(\star,\infty)}$ such that\footnote{One can also give bounds on the orders of $Q$ and $R$, but we will not need them.}
\[
 S=Q\circ_{\delta_\star} P+R \ \text{ with }  R=\sum_{i=0}^{p-1}R_j \circ_{\delta_\star} u_{\star}^j \ ,\  R_j\in 
K_{\eta^{-1}}.
\]
%
\end{Theorem}

\begin{Remark}
In the theorem above, consider for example the case $\star=0$. Following              \cite{Du2}*{Th\'eor\`eme 1.c.2}, the first step is to prove division in the case $P=u^p$. 
It is almost immediate to prove that one has right division, namely that for any $S\in\MM^{(0,\infty)}$ there are $Q,R$ such that $S=u^p\,Q+R$ and $R=\sum_{i=0}^{p-1} u_{\star}^jR_j$. 
 To prove left division as in the Theorem,
 notice that for $a(u)\in\C[[u]]$, $k\in\mathbb Z$ one has the formula
 \begin{equation}
 a(u)\,\eta^k= \sum_{i\geqslant 0}  (-1)^i{\binom{k}{i}}\cdot \eta^{k-i}\circ_{\delta_0} \delta_{0}^i(a(u)) \,,
 \end{equation}
 where for all $k,i\in\mathbb Z$\ \footnote{In fact, formula (4) is valid for any derivation $\delta$. One has also
 \[\eta^{k} \circ_{\delta} a(u)=\sum_{i\geqslant 0}{\binom{k}{i}}\,
 \delta^{i}(a(u)) \,\eta^{k-i}\,.\]}
 \[
  {\binom{k}{i}}=\frac{k(k-1) \cdots(k-i+1)}{i !}.
 \]
If $a(u)= u^{\ell}$, then $\delta_{0}^i(u^{\ell})=\ell (\ell-1)\cdots (\ell-i+1) u^{\ell}$,
 and so from the quotient and the remainder of right division we obtain those for left division. 
Similar considerations apply to the rest of the proof. 
\end{Remark}

\begin{Definition}
\begin{enumerate}
 \item  $\MM^{(s,\infty)}$: Let $\mathfrak H$ denote the category of holonomic $\C[[x]]\langle\partial_x\rangle$-modules, $\D$ the category of difference $K_{\theta}$-modules. For any  
$s\in \C\smallsetminus\{0\}$, we define a functor
\[
 \mathfrak M^{(s,\infty)}(\bullet): \mathfrak H \longrightarrow \D
\]
as follows: 
Consider the ring homomorphism
$\varphi_s:\C[[x]]\langle \partial_x \rangle \longrightarrow  \mathfrak M^{(s, \infty)}$
defined by $x \longmapsto u$ and  $(x+s)\partial_x \longmapsto -\eta$. Then, 
$\mathfrak M^{(s, \infty)}$ is a right $\C[[x]]\langle \partial_x \rangle$-module
via the product
\begin{eqnarray*}
 \mathfrak M^{(s, \infty)} \times \C[[x]]\langle \partial_x \rangle &\longrightarrow & \mathfrak M^{(s, \infty)} \\
 \!\!\! (\,P\,,\, q\,) &\longmapsto & P\circ_{\delta_s} \varphi_s(q)
\end{eqnarray*}

If $M$ is a $\C[[x]]\langle\partial_x\rangle$-module, we put $\mathfrak M^{(s,\infty)}(M)= \mathfrak M^{(s, \infty)}\otimes_{\C[[x]]\langle \partial_x \rangle} M $. If we set $\theta=\eta^{-1}$, then $\mathfrak M^{(s, \infty)}(M)$ is a vector space over $K_\theta$, where scalar multiplication is given by $f(\theta) \cdot (P(u,\eta)\otimes m)= (f(\eta^{-1})\circ_{\delta_s}P(u,\eta))\otimes m$. 

It is not difficult to check that the map
$\Phi_s: \mathfrak M^{(s, \infty)}(M) \longrightarrow \mathfrak M^{(s, \infty)}(M)$ defined by $ P(u,\eta) \otimes m \longmapsto (u+s)\cdot P(u,\eta)\otimes m$  endows $\mathfrak M^{(s, \infty)}(M)$ with a structure of $K_{\theta}$-difference module. It will follow from Theorem \ref{canex} and Lemma \ref{l3} below that this $K_{\theta}$-vector space has finite dimension (see Remark \ref{fd}), and so $\MM^{(s,\infty)}(\bullet)$ is a well-defined functor.

 \item  $\MM^{(0,\infty)}$:  
If $\mathfrak H'$ denote the category of  formal connections, we define a functor 
\[
  \mathfrak M^{(0,\infty)}(\bullet): \mathfrak H' \longrightarrow \D
\]
as follows: As before, $\mathfrak M^{(0, \infty)}$ is a right $\C[[x]]\langle x\partial_x \rangle$-module
via the ring homomorphism
$\varphi_0:\C[[x]]\langle x\partial_x \rangle \longrightarrow  \mathfrak M^{(0, \infty)}$
defined by $x \longmapsto u$ and  $x\partial_x \longmapsto -\eta$.

Let $M$ be a formal connection and
consider its canonical good filtration $V_{\ast}M$, recall that $V_0M$ is a $V_0\C[[x]]\langle \partial_x\rangle=\C[[x]]\langle x\partial_x\rangle$-module. In the tensor product
$\mathfrak M^{(0, \infty)}\otimes_{\C[[x]]\langle x\partial_x \rangle} V_0M$ we have a $K_{\theta}$--vector space structure, defined as in the previous case, and a difference operator $\Phi_0$ given by 
\[
 Q(u,\eta)\otimes m \longmapsto  u \cdot Q(u,\eta)\otimes  m \,.
\]
However, $\Phi_0$ might not be invertible, and then we would not have a difference module
as defined before. So, we invert it formally and set
\[
 \MM^{(0,\infty)}(M):=\mathfrak M^{(0, \infty)}\otimes_{\C[[x]]\langle x\partial_x \rangle} V_0M\otimes_{K_{\theta}\langle \Phi_0\rangle}K_{\theta}\langle \Phi_0,\Phi_0^{-1}\rangle .
\]
As in i), this is a finitely dimensional $K_{\theta}$-vector space (see Remark \ref{fd} below).

 \item  $\MM^{(\infty,\infty)}$: Finally, we define 
\[
   \mathfrak M^{(\infty,\infty)}(\bullet): \mathfrak H' \longrightarrow \D
\]
as in ii), replacing $u$ by $u^{-1}$. Namely, $\mathfrak M^{(\infty, \infty)}$ is a right $\C[[x]]\langle x\partial_x \rangle$-module
via the ring homomorphism
$\varphi_\infty:\C[[x]]\langle x\partial_x \rangle \longrightarrow  \mathfrak M^{(\infty, \infty)}$
defined by $x \longmapsto u^{-1}$ and  $x\partial_x \longmapsto -\eta$. Given a  formal connection $M$, 
we consider in the $K_{\theta}$--vector space $\mathfrak M^{(\infty, \infty)}\otimes_{\C[[x]]\langle x\partial_x \rangle} V_0M$ the difference operator $\Phi_{\infty}$ given by 
\[
 Q(u,\eta)\otimes m \longmapsto u^{-1}\cdot Q(u,\eta)\otimes  m\,.
\]
and we put $\mathfrak M^{(\infty,\infty)}(M):=\mathfrak M^{(\infty, \infty)}\otimes_{\C[[x]]\langle x\partial_x \rangle} V_0M\otimes_{K_{\theta}\langle \Phi_\infty\rangle}K_{\theta}\langle \Phi_\infty,\Phi_\infty^{-1}\rangle $.  Finitely dimensionality is proved as in the previous cases.

\end{enumerate}
The functors $ \mathfrak M^{(s,\infty)}(\bullet)$ send a morphism $f: M \longrightarrow N$ to $Id\otimes f$. The functors $\mathfrak M^{(0,\infty)}(\bullet)$ and  $\mathfrak M^{(\infty,\infty)}(\bullet)$, send $f$ to $Id\otimes f\otimes Id$, this is well defined by functoriality of the canonical good filtration.
\end{Definition}

\begin{Proposition} \label{flat}
 The ring homomorphisms $\varphi_{s}$ ($s\in\C$) and $\varphi_{\infty}$ are flat.
\end{Proposition}
\begin{Proof}
It is proved similarly as in \cite{Bj}*{Chap. 5, \S 5}, considering the $\eta$-order filtration in the spaces $\MM^{(\star,\infty)}$,  $\star\in\C\cup\{\infty\}$. $\Box$
\end{Proof}

\begin{Remark}
It follows immediately from this Proposition and from \cite{Sab}*{Corollary 6.1.3} that the functors  $\MM^{(\star,\infty)}(\bullet)$  are exact. 
\end{Remark}

We want to study  formal slopes and horizontal zeros of  local Mellin transforms. For modules with punctual support, we have
\begin{Lemma} \label{punk}
 Let $M$ be a holonomic $\C[[x]]\langle \partial_x \rangle$-module supported only at zero.
 If $s\in \C\smallsetminus\{0\}$, then $\mathfrak M^{(s,\infty)}(M)$ is purely of slope zero, $Hor(\MM^{\,(s,\infty)} (M))\subset \{s\}$ and $\dim_{K_{\theta}} \MM^{\,(s,\infty)} (M)= \mu(M)$.
 
\end{Lemma}
\begin{Proof}
We can assume that $M= \C[[x]]\langle \partial_x \rangle/\C[[x]]\langle \partial_x \rangle\cdot x^m$ for some $m\geqslant 1$, then $\mu(M)=m$. 
We have isomorphisms of difference modules
\[
 \mathfrak M^{(s,\infty)}(M)=\frac{\mathfrak M^{(s,\infty)}}{\mathfrak M^{(s,\infty)}\cdot u^m}\cong \frac{K_{\theta}\langle \Phi\rangle}{K_{\theta}\langle \Phi\rangle (\Phi-s)^m}\,,
\]
and the assertions follow easily. $\Box$
\end{Proof}

To treat the general case, we will need the following theorem, proved in \cite{Mal}*{(5.1)} (see also \cite{Ka}*{Theorem 2.4.10}):

\begin{Theorem} \label{canex}
Let $N$ be a holonomic $\mathbb C[[z]]\langle \partial_z\rangle$-module. Then, there exists a holonomic
$\mathbb C[z]\langle \partial_z\rangle$-module $\mathbb M_{N}$ with no singularity on $\mathbb A^ 1_{\mathbb C}\smallsetminus \{0\}$, regular at infinity and such that $(\mathbb M_N)_0 \cong N$.
\end{Theorem}

If $s\in \mathbb C$, $\tau_s: \mathbb A^1_{\mathbb C} \longrightarrow \mathbb A^1_{\mathbb C}$ is the translation $z\mapsto z-s$
and $i: \mathbb A_{\mathbb C}^ 1\smallsetminus \{0\} \hookrightarrow \mathbb A^1_{\mathbb C}$ is the inclusion map, we put 
$\mathbb M_{N}^s=i^{\ast}\tau_s^{\ast}(\mathbb M_N)$. If $j: \mathbb A_{\mathbb C}^ 1\smallsetminus \{0\} \longrightarrow \mathbb A_{\mathbb C}^ 1\smallsetminus \{0\}$ is the inversion
$z\mapsto z^{-1}$, we put $\mathbb M^{\infty}_N=j^{\ast}i^{\ast}(\mathbb M_N)$,
where $\tau_s^{\ast}, i^{\ast}, j^{\ast}$ denote inverse images as $\mathcal D$-modules.

\begin{Lemma} \label{l2}
\begin{enumerate}
 \item The slopes of  $\MM(\mathbb M^0_N)_{\infty}$, the germ at infinity of the global Mellin transform of $\mathbb M^0_N$, are strictly negative.
 \item For $s\in \mathbb C\smallsetminus \{0\}$, $\MM(\mathbb M_N^s)_{\infty}$ has only slope $0$ and $Hor(\MM(\mathbb M^s_N)_{\infty})\subset\{s\}$.
 \item The slopes of  $\MM(\mathbb M^{\infty}_N)_{\infty}$ are strictly positive.
\end{enumerate}
\end{Lemma}
\begin{Proof} Since $\M_N^0$ is regular at infinity, the Newton polygon of $\M_N^0$ has only sides of non-negative slope (corresponding to the slopes of $N$) and no vertical side (since $\M_N^0$ has no singular points in $\mathbb A^1_{\mathbb C}\smallsetminus \{0\}$). Then, item i) follows from Lemma \ref{pol}. The claims about slopes in ii) and iii) is proved in the same way.

In case ii), consider an exact sequence 
\[
 0 \longrightarrow \mathbb K \longrightarrow \M'\longrightarrow \M^s_N\longrightarrow 0\,
\]
as in the proof of Lemma \ref{pol}.
Since the only singular point of $\M^s_N$ is $s\in \mathbb C\smallsetminus \{0\}$ and $\M^s_N$ is of exponential type in the sense of \cite{Mal}*{Chapitre XII}, $\M'$ will be of exponential type as well, and if $\mathbb K=\mathbb C[z,z^{-1}]\langle z\partial_z\rangle/ \mathbb C[z,z^{-1}]\langle z\partial_z\rangle\cdot q(z)$ with $q(z)\in \mathbb C[z]$,  we will have 
\[
 P=(z-s)^{m}q(z)\partial_z^n+\dots \ \ \text{,\ where }\ m+\deg(q)=\deg_z(P) , \  n=\deg_{\partial_z}(P)\,,
\]
where the coefficients of $\partial^i_z$ are in $\mathbb C[z]$ for all $i\geqslant 0$. 
A computation using just the definitions and (\ref{2}) shows that $Hor(\MM(\M_N^s)_{\infty})\subset\{s\}$. $\Box$
\end{Proof}

\begin{Lemma} \label{l3}
Let $\mathbb M$ be a holonomic $\C[z,z^{-1}]\langle z\partial_z\rangle$-module. 
If $\star\in \C\cup\{\infty\}$, consider the maps $\Xi_\star:\ \mathfrak M (\mathbb M)_{\infty} \longrightarrow \mathfrak M^{(\star,\infty)}(\mathbb M_\star)$ defined as follows

  \begin{enumerate}
   \item If $\star=s\in\C\smallsetminus\{0\}$, then
   \begin{eqnarray*}
 \!\!\!\!\!\!\!\!\!\!\Xi_s:\ \mathfrak M (\mathbb M)_{\infty} = K_{\theta}\langle \Phi,\Phi^{-1}\rangle\otimes_{\C[z,z^{-1}]\langle z\partial_z\rangle} \M  & \longrightarrow &
 \mathfrak M^{(s,\infty)}(\mathbb M_s)\\ \!\!\!
 \Phi^k \otimes m \  & \longrightarrow & (u+s)^k  \otimes m  \,.
   \end{eqnarray*}
   \item  If $\star=0$,
   \begin{eqnarray*}
 \!\!\!\!\!\!\!\!\!\!\Xi_0:\ \mathfrak M (\mathbb M)_{\infty} = K_{\theta}\langle \Phi, \Phi^{-1}\rangle\otimes_{\C[z]\langle z\partial_z\rangle} V_0\M  & \longrightarrow &
 \mathfrak M^{(0,\infty)}(\M_0)\\
  \Phi^k \otimes m \  & \longrightarrow &   1\otimes m  \otimes \Phi_0^k\,.
\end{eqnarray*}
\item If $\star=\infty$,
\begin{eqnarray*}
 \!\!\!\!\!\!\!\!\!\!
 \Xi_\infty:\ \mathfrak M (\mathbb M)_{\infty} = K_{\theta}\langle \Phi, \Phi^{-1}\rangle\otimes_{\C[z^{-1}]\langle z^{-1}\partial_{z^{-1}}\rangle} \widetilde{V_0\M}  & \longrightarrow &
 \mathfrak M^{(\infty,\infty)}(\M_\infty)\\
  \Phi^k \otimes m \  & \longrightarrow &    1\otimes m \,\otimes\Phi_{\infty}^k.
\end{eqnarray*}
 \end{enumerate}
In all three cases we extend by $K_{\theta}$-linearity. 
 Then, the maps $\Xi_{\star}$ are epimorphisms of $K_{\theta}$--difference modules, functorial on $\M$.
\end{Lemma}

\begin{Proof}
That the given maps are morphisms of $K_{\theta}$-difference modules, functorial on $\M$, follows from the definitions. To see they are onto,  take an epimorphism of holonomic $\C[z,z^{-1}]\langle z\partial_z\rangle$-modules
\[
\frac{\C[z,z^{-1}]\langle z\partial_z\rangle}{\C[z,z^{-1}]\langle z\partial_z\rangle \cdot P(z,z\partial_z)} \longrightarrow \M \longrightarrow 0
\]
Since both local and global Mellin transforms are exact functors, by functoriality  we can assume that 
\[
\mathbb M = \frac{\C[z,z^{-1}]\langle z\partial_z\rangle}{\C[z,z^{-1}]\langle z\partial_z\rangle \cdot P(z,z\partial_z)} .
\]

{\it \underline{Case} $\Xi_s,\ s\in\C\smallsetminus \{0\}$}: We have
\[
 \mathfrak M (\mathbb M)_{\infty} = \frac{K_{\theta} \langle \Phi,\Phi^{-1}\rangle}{K_{\theta} \langle \Phi,\Phi^{-1}\rangle\cdot P(\Phi, -\theta^{-1})}
\]
and
\[
\MM^{(s,\infty)}(\mathbb M_{s}) =\frac{\MM^{(s,\infty)}}{\MM^{(s,\infty)}\cdot P(u+s, -\eta)}\,.
\]
The map $\Xi_{s}$ is given by $[a(\theta)\,\Phi^k]\longmapsto [ a(\eta^ {-1})\circ_{\delta_s} (u+s)^k]$. Surjectivity follows from the division theorem.

{\it \underline{Case} $\Xi_0$}: 
Let $b(s)\in\C[s]$ be the Bernstein polynomial of the canonical good $z$-filtration $V_{\ast}\M$ of $\M$ (see for example \cite{Sab}*{I, Section 6}). By definition of $b(s)$, there is a $P'\in V_{-1}\C[z,z^{-1}]\langle z\partial_z\rangle$ such that $Q(z,z\partial_z)=b(z\partial_z)-zP'$ annihilates $\M$ . So, we have an epimorphism
\[
 \frac{\C[z,z^{-1}]\langle z\partial_z\rangle}{\C[z,z^{-1}]\langle z\partial_z\rangle \cdot Q(z,z\partial_z)} \longrightarrow  \M \longrightarrow 0
\]
and, as before, we can assume
\[
 \M= \frac{\C[z,z^{-1}]\langle z\partial_z\rangle}{\C[z,z^{-1}]\langle z\partial_z\rangle \cdot Q(z,z\partial_z)} .
\]
In this case, and by the special shape of $Q(z,z\partial_z)$ (see \cite{Sab}*{Exercise 6.1.4 (3)}), the canonical good $z$-filtration of $\M$ is given by
\[
 V_k\M= \frac{V_k\C[z,z^{-1}]\langle z\partial_z\rangle}{V_k \C[z,z^{-1}]\langle z\partial_z\rangle \cdot Q}.
\]
Then we have 
\begin{eqnarray*}
\MM^{(0,\infty)}(\mathbb M_{0}) & \cong &  \frac{\MM^{(0,\infty)}}{\MM^{(0,\infty)}\cdot Q(u, -\eta)}\otimes_{K_{\theta}\langle \Phi_0\rangle} K_{\theta}\langle \Phi_0,\Phi_0^{-1}\rangle \,,\\ \ \\
\MM(\mathbb M)_{\infty}&\cong &\frac{K_\theta\langle\Phi,\Phi^{-1}\rangle}{K_\theta \langle\Phi,\Phi^{-1}\rangle\cdot Q(\Phi,-\theta^{-1})}\,,
\end{eqnarray*}
and, via this isomorphisms, the map $\Xi_0$ is given by $[a(\theta)\,\Phi^k]\longmapsto  [ a(\eta^ {-1})] \otimes \Phi_0^k$. Again, one gets surjectivity from the division theorem.

{\it \underline{Case} $\Xi_\infty$}: It is analogous to the previous one, considering the canonical $z^{-1}$-filtration.
$\Box$
\end{Proof}

\begin{Remark} \label{fd}
 If $M$ is a formal connection and $\star\in\{0,\infty\}$, it follows from this lemma applied to $\M^\star_M$ that  the $K_{\theta}$-vector spaces $\MM^{(\star,\infty)}(M)$ 
 are finite dimensional. By the same argument, if $M$ is a holonomic $\C[[x]]\langle\partial_x\rangle$-module and $s\in\C\smallsetminus\{0\}$, then $\MM^{(s,\infty)}(M)$
 is finite dimensional as well.
\end{Remark}

\begin{Proposition} \label{locme}

\begin{enumerate}
 \item   Let $M$ be a formal connection. Then: 
 \begin{enumerate}
\item  All slopes of $\MM^{\,(0,\infty)}(M)$ are strictly negative and \\ $\dim_{K_{\theta}} \MM^{\,(0,\infty)} (M)\geqslant  \mathrm{irr}(M)$. 
 \item  All slopes of $\MM^{\,(\infty,\infty)}(M)$ are strictly positive and \\ $\dim_{K_{\theta}} \MM^{\,(\infty,\infty)} (M)\geqslant \mathrm{irr}(M)$.
 \end{enumerate}
 \item Let $M$ be a holonomic $\C[[x]]\langle \partial_x \rangle$-module. If $s\in \C\smallsetminus\{0\}$, then $\MM^{\,(s,\infty)}(M)$ has only slope zero,
   $Hor(\MM^{\,(s,\infty)} (M))\subset \{s\}$, and  $\dim_{K_{\theta}} \MM^{\,(s,\infty)} (M)= \mu(M)$. 

\end{enumerate}
\end{Proposition}
\begin{Proof} 
All assertions about slopes and the assertion in (ii) about horizontal zeros  follow from Lemma \ref{l2} and Lemma \ref{l3} applied to $\mathbb M^{\star}_M$.  We prove the remaining claim in (a): We can assume that
\[
M = \frac{\mathbb C[[x]]\langle \partial_x\rangle}{\mathbb C[[x]]\langle \partial_x\rangle \cdot P(x,x\partial_x)},
\]
where
\[
 P(x,x\partial_x)=a_d(x)(x\partial_x)^d + \dots + a_0(x) \,,
\]
$a_i(x)\in \mathbb C[[x]]$ for all $i\geqslant 0$, at least one of the $a_i(x)$ is a unit, and 
$a_d(x)\neq 0$. 
Then, the irregularity of $M$ is $\text{ord}_x (a_d(x))$. 

We have a good filtration on $M$ defined  by
\[
 U_k M=\frac{ V_k\mathbb C[[x]]\langle \partial_x\rangle}{V_k \mathbb C[[x]]\langle \partial_x\rangle\cdot P}  \ \ \text{ for } k\in\mathbb Z
\]
In general, this will not be the canonical good filtration $V_{\ast}M$ because the roots of the Bernstein polynomial of $U_{\ast}M$ do not need to be contained in $\Sigma$. 
However, we will have $U_{-k}M\subset V_0M$, for some $k\geqslant 0$ (see e.g. \cite{Sab}*{Exercise 5.13, 5}). Also, there is an isomorphism of $V_0\mathbb C[[x]]\langle \partial_x\rangle$-modules
\[
 U_{-k}M= \frac{V_{-k} \C[[x]]\langle \partial_x\rangle}{V_{-k} \C[[x]]\langle \partial_x\rangle\cdot P} = \frac{V_{0} \C[[x]]\langle \partial_x\rangle\, x^k}{V_{0} \C[[x]]\langle \partial_x\rangle\, x^k \cdot P}\cong \frac{V_{0} \C[[x]]\langle \partial_x\rangle}{V_{0} \C[[x]]\langle \partial_x\rangle\, x^k \cdot P \cdot x^{-k}}\,,
\]
where the second equality follows from $V_{-k} \C[[x]]\langle \partial_x\rangle=V_{0} \C[[x]]\langle \partial_x\rangle\, x^k $ and the last isomorphism holds because right multiplication by $x^k$ is a morphism of left $V_0\mathbb \C[[x]]\langle \partial_x\rangle$-modules.
If we set $Q= x^k \cdot P \cdot x^{-k}$, then one sees that $Q\in V_0\C[[x]]\langle \partial_x\rangle$ and, perhaps increasing $k$, we have
\[
 Q(x,x\partial_x)= b_d(x)(x\partial_x)^d + \dots + b_0(x) \,,
\]
where $b_i(x)\in \C [[x]]$ for all $i\geqslant 0$, at least one of the $b_i(x)$ is a unit, and $b_d(x)=a_d(x)$.\footnote{All this follows easily from the identity 
$
x^k(x\partial_x)^dx^{-k}=(x\partial_x-k)^d\ k,d\in\mathbb N
$
To insure that at least one of the $b_i(x)$ is a unit, it might be necessary to increase the value of $k$. } We put 
\begin{equation} \label{mq}
 \MM_Q:= \MM^{(0,\infty)}\otimes_{\C[[x]]\langle x\partial_x \rangle} \frac{V_{0} \C[[x]]\langle \partial_x\rangle}{V_{0} \C[[x]]\langle \partial_x\rangle\, Q}\cong \frac{ \MM^{(0,\infty)}}{ \MM^{(0,\infty)} \cdot Q(u,-\eta)}.
\end{equation}

By flatness of $\MM^{(0,\infty)}$ as a $\C[[x]]\langle x\partial_x\rangle$-module we have
an injective morphism
 $\MM_Q \hookrightarrow  \MM^{(0,\infty)}\otimes V_0M$,
and it follows from the division theorem that  
\[
\dim_{K_{\theta}} \MM_Q
=\text{ord}_x (b_d(x))=\mathrm{irr}(M)\,.
\] 
Thus, to prove the desired inequality it suffices to show that
\begin{equation} \label{ine}
 \dim_{K_{\theta}} \MM_Q  \geqslant \dim_{K_{\theta}}  \MM_Q\otimes_{K_{\theta}\langle \Phi_0 \rangle} K_{\theta}\langle \Phi_0,\Phi_0^{-1}\rangle .
\end{equation}
We claim that the action of $\Phi_0$ on $\MM_Q$ is bijective: 
Since $\MM_Q$ is finite dimensional over $K_{\theta}$
and $\phi: K_{\theta} \longrightarrow K_{\theta}$ is an automorphism,  it is enough to prove that $\Phi_0$ is injective. Via the isomorphism (\ref{mq}), the action of $\Phi_0$ is
given by $[R(u,\eta)] \longmapsto [u\,R(u,\eta)]$. If $u\,R= T\,Q$ for some $T\in\MM^{(0,\infty)}$, then, dividing $T$ by $u$ we have $u\,R= (uT'+T'')\,Q$ and, if we take classes in the quotient of $\MM^{(0,\infty)}$ by the right ideal generated by $u$, we obtain an equality 
$
0=\widetilde{T''}\cdot \widetilde{Q} 
$
in the field $K_{\theta}=K_{\eta^{-1}}$. But $\widetilde{Q}\neq 0$ because at least one of the $b_i(u)$ is a unit in $\C[[u]]$, so $T''=0$ and it follows that $R\in\MM^{(0,\infty)}Q$, as desired.

Finally, since $\Phi_0: \MM_Q \longrightarrow \MM_Q$ is bijective, 
$\MM_Q \cong \MM_Q\otimes_{K_{\theta}\langle \Phi_0 \rangle} K_{\theta}\langle \Phi_0,\Phi_0^{-1}\rangle $
and (\ref{ine}) follows. The proof of the inequality in item (b) is analogous.

We prove next the last assertion in (ii): There is an exact sequence of $\C [[x]]\langle \partial_x \rangle$-modules
\[
 0 \longrightarrow K_1\longrightarrow M \longrightarrow M'=\frac{\C [[x]]\langle \partial_x\rangle}{\mathbb C[[x]]\langle \partial_x\rangle \cdot P(x,x\partial_x)} \longrightarrow K_2 \longrightarrow 0
\]
where $K_1,K_2$ are supported at zero and $M'$ is a formal connection.  To prove that $\dim_{K_{\theta}} \MM^{\,(s,\infty)} (M)= \mu(M)$, it suffices to prove this same equality for $K_1,K_2$ and $M'$. For $K_1$ and $K_2$, see Lemma \ref{punk}, for $M'$ we have
\[
 \MM^{\,(s,\infty)}(M')= \frac{ \MM^{\,(s,\infty)}}{ \MM^{\,(s,\infty)}\cdot P\left(u, \left(\sum_{i\geqslant 0}(\frac{-u}{s})^ {i+1}\right)\eta\right)},
\]
and $P(u, \left(\sum_{i\geqslant 0}(\frac{-u}{s})^ {i+1}\right)\eta)=c_d(u)\eta^ d + c_{d-1}(u) \eta^ {d-1}+ \dots$ with $\text{ord}_u(c_d(u))=d+\text{ord}_u(a_d(u))$. Again by the division theorem, we have 
$
\dim_{K_{\theta}} \MM^{\,(s,\infty)}(M')=d+\text{ord}_u(a_d(u))=\mu(M')
$, as desired. 
%
$\Box$
\end{Proof}


{\bf 4. Formal stationary phase for the Mellin transform.}

In this section we prove the main theorem stated in the introduction:

\begin{Theorem}
 Let $\mathbb M$ be a holonomic $\C[z,z^{-1}]\langle z\partial_z\rangle$-module. 
The map 
\[
\Xi=\!\!\! \!\!\!\!\!\bigoplus_{\star\in S(\mathbb M)\cup\{0,\infty\}}\, \!\!\! \!\!\!\!\!\! \Xi_{\star}: \,\MM(\M)_{\infty} \longrightarrow  \!\!\! \!\!\!\!\bigoplus_{\star\in S(\mathbb M)\cup\{0,\infty\}}  \!\!\! \!\!\!\!\MM^{(\star, \infty)}(\M_{\star})
\]
is a functorial isomorphism of difference modules over $K_{\theta}$. 
\end{Theorem}
\begin{Proof}

We show first that the map $\Xi$ is onto:
We can decompose $\MM(\M)_{\infty}$ according to slopes
\[
  \MM(\M)_{\infty} \cong  \MM(\M)_{\infty}^{>0} \oplus \MM(\M)_{\infty}^{=0}  \oplus \MM(\M)_{\infty}^ {<0} 
\]
and then, by Corollary  \ref{coslo}, Lemma \ref{l3} and Proposition \ref{locme}, it follows that it suffices to prove that 
the map
\[
\oplus_{s\in S(\mathbb M)}\, \Xi_s: \MM(\M)_{\infty}^{=0} \longrightarrow \bigoplus_{s\in 
S(\mathbb M)}\MM^{\,(s, \infty)}(\M_s)
\]
is onto. By the classification theorem, after a cyclic extension of $K_{\theta}$ we can assume that $ \MM(\M)_{\infty}^{=0}$ can be decomposed according to its horizontal zeros. But from Proposition \ref{locme} we have 
$Hor(\MM^{\,(s, \infty)}(\M_s))\subset \{s\}$ for $s\in\mathbb S(\mathbb M)$ and then by Lemmas \ref{hor} and \ref{l3}
the surjectivity of $\oplus_{s\in S(\mathbb M)}\, \Xi_s$ follows.

Thus, it is enough to show that the dimension over $K_{\theta}$ of the source of $\Xi$ is smaller or equal than the dimension of its target: By Lemma \ref{pol}, the height of 
the Newton polygon of $\M$ coincides with the width of the Newton polygon of $\MM(\M)_{\infty}$. So, we have\footnote{To my knowledge, this formula was first proved by C. Sabbah, using a different method (unpublished, 
but see \cite{LS2} for a special case).}
\[
 \dim_{K_{\theta}}(\MM(\M)_{\infty})=  \mathrm{irr}(\mathbb M_0) + \mathrm{irr}(\mathbb M_{\infty}) + \sum_{s\in S(\mathbb 
M)}\mu(\mathbb M_s)\,
\]
and then, by Proposition \ref{locme} we are done. 
$\Box$
\end{Proof}

Applying the theorem to modules of type $\M^0_N$, $\M^\infty_N$, we get

\begin{Corollary}
 The inequalities in Proposition \ref{locme}, items (a) and (b), are in fact equalities.
\end{Corollary}

\begin{Remark} 
Local Fourier transforms can be defined at the analytic level (this is well-known, take first an extension to a holonomic module on the affine $\mathbb C$-line, apply the global Fourier transform and localize at infinity), denote $\Phi^{(\star,\infty)}$ the functors so defined ($\star\in \mathbb C \cup \{\infty\}$). For differential modules of rank one in one variable, the formal and the analytic classification coincide, so  the formal stationary phase isomorphism for the Fourier transform \cite{GL}*{section 1} implies that if $\M$ is a holonomic module over the affine line, then there is an analytic isomorphism\footnote{See \cite{Mo}*{(5.11)} for a description, in terms of local data, of the determinant of the global Fourier transform (and not just of its germ at infinity).}
\begin{equation*}
 \det \mathfrak{Four}(\M)_{\infty} \cong \otimes_{\star\in\mathbb C \cup \{\infty\}} \det \Phi^{(\star, \infty)}(\M_\star)\ 
\end{equation*}
Following the analogous procedure, local Mellin transforms can also be defined at the analytic level. However, for difference modules of rank one, the analytic classification is much finer than the formal one (see \cite{vPS}*{10.2}), and therefore Theorem 4 above does not allow to derive a similar conclusion as in the Fourier case. If $\M$ is a module with regular singularities, then the analytic type of the determinant of its Mellin transform was determined in \cite{LS} (as explained in loc. cit., in fact only regularity at zero and at infinity is needed).

One could ask about possible $\ell$--adic analogues of the local Mellin transforms. The global Mellin transform does have a $\ell$-adic analogue, see \cite{GLoe}, and its determinant was computed in \cite{Loe}. Also, since in our approach the local Mellin transforms are a kind of slightly modified microlocalizations and since there is a good theory of $p$-adic microdifferential operators and local Fourier transforms (see \cite{Abe}, \cite{AbeMar}), one could also hope 
for a $p$-adic theory of local Mellin transforms, which might be related to the results in \cite{Loe2}.
\end{Remark}

%

\end{document}